\pdfoutput=1
\RequirePackage{ifpdf}
\ifpdf 
\documentclass[pdftex]{sigma}
\else
\documentclass{sigma}
\fi

\numberwithin{equation}{section}

\newtheorem{Theorem}{Theorem}[section]
\newtheorem*{Claim}{Claim}

\newtheorem{Lemma}[Theorem]{Lemma}
\newtheorem{Proposition}[Theorem]{Proposition}
\newtheorem{Setting}[Theorem]{Setting}
\newtheorem{defprop}[Theorem]{Definition/Proposition}

 { \theoremstyle{definition}
\newtheorem{Definition}[Theorem]{Definition}

\newtheorem{Remark}[Theorem]{Remark} }

\usepackage{tikz}

\DeclareMathOperator{\Ran}{Ran}
\DeclareMathOperator{\Ker}{Ker}
\DeclareMathOperator{\Coker}{Coker}
\DeclareMathOperator{\ind}{ind}
\DeclareMathOperator{\spann}{span}
\DeclareMathOperator{\Dim}{Dim}

\begin{document}

\allowdisplaybreaks

\newcommand{\arXivNumber}{2408.16696}

\renewcommand{\PaperNumber}{088}

\FirstPageHeading

\ShortArticleName{Perturbations of APS Boundary Conditions for Lorentzian Dirac Operators}

\ArticleName{Perturbations of APS Boundary Conditions\\ for Lorentzian Dirac Operators}

\Author{Lennart RONGE}

\AuthorNameForHeading{L.~Ronge}

\Address{Universit\"at Potsdam, Haus 9, Karl-Liebknecht-Str. 24-25, 14476 Potsdam, Germany}
\Email{\href{mailto:ronge@uni-potsdam.de}{ronge@uni-potsdam.de}}

\ArticleDates{Received May 20, 2025, in final form October 09, 2025; Published online October 20, 2025}

\Abstract{We study how far APS boundary conditions for a Lorentzian Dirac operator may be perturbed without destroying Fredholmness of the Dirac operator. This is done by developing criteria under which the perturbation of a compact pair of projections is a~Fredholm pair.}


\Classification{58J20; 58J45; 58J32; 47A53}

	\section{Introduction}
For Riemannian manifolds with boundary, the question which boundary conditions make a~Dirac operator Fredholm is well understood. However, in the Lorentzian case much less is known (see~\cite{BH} for what is known there). In general,  in order to find out whether or not a Dirac operator on a globally hyperbolic spacetime with spacelike boundary is Fredholm with given boundary conditions, one needs knowledge of the associated evolution operator, which is fairly hard to obtain.
	
For specific boundary conditions known as Atiyah--Patodi--Singer (APS) boundary conditions, the Lorentzian Dirac operator is always Fredholm and an index formula analogous to the Riemannian APS index theorem holds (see \cite{BSglobal,BSlocal}). In fact, with these boundary conditions, the operator is, loosely speaking, ``as Fredholm as possible'' and this makes it possible to extend the result to boundary conditions that are close to APS boundary conditions. In \cite[Proposition~4.5]{BH}, it is shown that boundary conditions given by the graphs of certain operators (such that the zero operators correspond exactly to APS-conditions) still make the Dirac operator Fredholm as long as the product of the operators' norms is less than some unspecified constant~${\varepsilon>0}$.
	
The goal of this note is to improve this condition and study in more detail, how far boundary conditions may be allowed to deviate from APS conditions without becoming non-Fredholm. We find that the $\varepsilon$ in \cite[Proposition~4.5]{BH} may be chosen to be 1 (Theorem~\ref{finalgraph}), which is the best constant that one could hope to achieve. We also provide a different condition in terms of the difference in the norms (modulo compact operators) of the boundary projections and the APS projections (Theorem~\ref{finalCnorm}).
	
	All the Lorentzian geometry we will need is already present in \cite{BH}. The actual work lies purely in the study of general pairs of projections, for which we derive Fredholmness criteria. These criteria, applied to the setting of \cite{BH} will then yield the improved results.
	
	The structure is as follows: We provide background information in Section~\ref{s2}, derive Fredholmness criteria for general pairs of projections in Section~\ref{s3} and use them in Section~\ref{s4} to obtain results for Lorentzian boundary value problems.
	
\section{Background}\label{s2}
	\subsection{Fredholm pairs and Calkin norm}
	For this subsection, let $H$ and $H'$ be complex Hilbert spaces. Let $\mathcal{B}(H)$ denote bounded operators on $H$, $\mathcal{K}(H)$ denote compact operators on $H$ and $\mathcal C(H):= \mathcal{B}(H)/\mathcal{K}(H)$ the Calkin algebra of~$H$. Let $\mathcal{B}(H,H')$ and $\mathcal{K}(H,H')$ denote bounded operators and compact operators from $H$ to~$H'$.
	We recall the definition and well-known characterizations of Fredholm operators.
	\begin{defprop}
		An operator $T\in \mathcal{B}(H,H')$ is called Fredholm, if one of the following equivalent conditions is satisfied:
		\begin{itemize}\itemsep=0pt
			\item $\Ker(T)$ and $\Coker(T):=H'/\Ran(T)$ are finite-dimensional.
			\item $\Ker(T)$ and $\Ker(T^*)$ are finite-dimensional and $T$ has closed range
		\end{itemize}
		If $T\in \mathcal{B}(H)$, $T$ is Fredholm if and only if its equivalence class in $\mathcal{C}(H)$ is invertible. The operator $T^*$ is Fredholm if and only if $T$ is Fredholm. The index of a Fredholm operator $T$ is defined by
		\[\ind(T):=\Dim\Ker(T)-\Dim\Coker(T)=\Dim\Ker(T)-\Dim\Ker(T^*).\]
	\end{defprop}
	A central role in this paper is played by Fredholm pairs of projections.
	\begin{Definition}
		Whenever we refer to a projection (in $H$) in this paper, we shall mean an orthogonal projection onto a closed subspace (of $H$).
	\end{Definition}
	\begin{Definition}
		A pair $(P,Q)$ of projections in $H$ is called Fredholm (of index $k$), if the restricted operator
		\[Q\colon\ \Ran(P)\rightarrow\Ran(Q)\]
		is Fredholm (of index $k$). We will denote this restricted operator by $Q_P$.
	\end{Definition}
	 For the source of this definition and further information on Fredholm pairs, we refer the reader to \cite{ASS}. The facts that we will need are the following (\cite[Theorem~3.4 and Proposition~3.1]{ASS} and \cite[Lemma~3.2]{LFP}).
	 \begin{Proposition}
		\label{Fredholmprops}
		Let $(P,Q)$ be a Fredholm pair of projections.
		\begin{enumerate}\itemsep=0pt
			\item[$(1)$] The pair $(Q,P)$ is Fredholm and $\ind(Q,P)=-\ind(P,Q)$.
			\item[$(2)$] For any unitary $U$, the pair $\bigl(UPU^{-1}, UQU^{-1}\bigr)$ is Fredholm with the same index as $(P,Q)$.
			\item[$(3)$] If $Q'$ is a projection and $Q'-Q$ is compact, then $(P,Q')$ is Fredholm and
			\[\ind\bigl(P,Q'\bigr)=\ind(P,Q)+\ind\bigl(Q,Q'\bigr).\]
			\item[$(4)$] The pair $(1-P,1-Q)$ is Fredholm and
			\[\ind(1-P,1-Q)=-\ind(P,Q).\]
			\item[$(5)$] If $((P(t),Q(t)))_{t\in [0,1]}$ is a continuous path of Fredholm pairs, then
			\[\ind(P(1),Q(1))=\ind(P(0),Q(0)).\]
    \end{enumerate}
		\end{Proposition}
	
The norm of the Calkin algebra will be important for our study of Fredholm pairs. It can be viewed in different ways.
	\begin{defprop}		\label {cnorm}
		For $T\in\mathcal{B}(H)$, denote by $\|T\|_C$ the norm of the equivalence class of $T$ in $\mathcal{C}(H)$. This is equal to each of the following:
		\begin{enumerate}\itemsep=0pt
			\item[$(1)$] 
            $\|T\|_1:=\inf\{\|S\| \mid S\in\mathcal{B}(H),\, S-T\in\mathcal{K}(H)\}$.
			\item[$(2)$] 
            $\|T\|_2:=\inf\{\|S\| \mid S\in\mathcal{B}(H),\, S-T\text{ has finite rank}\}$.
			\item[$(3)$] 
            $\|T\|_3:=\inf\{\|T|_V\|\mid V\subseteq H \text{ is a closed subspace of finite codimension}\}$.
		\end{enumerate}
	\end{defprop}
	\begin{proof}
		$\|T\|_1$ is the definition of the quotient norm.
		
		Proof of $\|T\|_1=\|T\|_2$:
				As $\|T\|_1$ is an infimum over a larger set, we immediately get
		\[\|T\|_1\leq\|T\|_2.\]
		Let $\varepsilon>0$ and let $K$ be a compact operator such that \[\|T+K\|<\|T\|_1+\varepsilon.\] As a compact operator, $K$ can be approximated up to $\varepsilon$ by a finite rank operator $K'$. We obtain
		\[\bigl\|T+K'\bigr\|<\|T\|_1+2\varepsilon.\]
		As the left-hand side is an upper bound for the infimum defining $\|T\|_2$, we obtain
		\[\|T\|_2\leq \|T\|_1+2\varepsilon.\]
		As $\varepsilon$ was arbitrary, the two infima coincide.
		
		Proof of $\|T\|_2=\|T\|_3$:
For $K$ finite rank, we have
        \[
        \|T+K\|\geq \bigl\|(T+K)|_{\Ker(K)}\bigr\|= \bigl\|T|_{\Ker(K)}\bigr\|.
        \]
        Taking the infimum over all finite rank $K$ and using that $\Ker(K)$ has finite codimension, we obtain $\|T\|_2\geq\|T\|_3$.
		
		Conversely, if $V$ has finite codimension, let $P_V$ denote the projection onto $V$. As $T-TP_V$ has a kernel with finite codimension,  it must have finite rank. As we have $\|T|_V\|=\|TP_V\|$, we can conclude (taking the infimum over all $V$) that  $\|T\|_3\geq\|T\|_2$.
	\end{proof}

	An alternative criterion for Fredholmness that will be important is the following (see \cite[Corollary~5.3.14]{DSW}, in a different form \cite[Proposition~3.2]{ASS}).
	\begin{Proposition}
		\label{FredDiff}
		A pair of projections $(P,Q)$ is a Fredholm pair if and only if
		$\|P-Q\|_C<1$.
	\end{Proposition}
	In particular, whether a pair is Fredholm (but not its index) is determined by the classes of the projections in the Calkin algebra. If the difference $P-Q$ is compact, then it has norm $0$ in the Calkin algebra, thus the pair is in some sense ``as Fredholm as possible''.
	\subsection{Lorentzian index theory}
	\label{Lorentzintro}
	We use the same setting as \cite{BH}. As we mostly talk about abstract Fredholm pairs anyways, we keep our description of the Lorentzian setting brief and refer the interested reader to \cite[Section~2]{BH} for more detailed explanations. More generally, any setting that has a unitary evolution operator, allows the definition of APS boundary conditions and  in which Theorems \ref{bdyindex} and \ref{APScomp} hold will work.

	\begin{Setting}
We assume the following Lorentzian setting:
	\begin{itemize}\itemsep=0pt
            \item $M$ is an even-dimensional spatially compact globally hyperbolic time-oriented Lorentzian spin manifold.
			\item $\Sigma_0$ and $\Sigma_1$ are smooth spacelike Cauchy hypersurfaces, with $\Sigma_0$ in the past of $\Sigma_1$.
			\item $M_0$ is the spacetime region between $\Sigma_0$ and $\Sigma_1$.
			\item $E$ is a Hermitian vector bundle over $M$ with a connection.
			\item $S(M)\otimes E$ is the associated twisted spinor bundle, $S^R(M)\otimes E$ is the subbundle of right-handed spinors.
			\item $D$ is the twisted Dirac operator acting on right-handed spinor fields of the twisted bundle.
			\item $Q\colon L^2\bigl(\Sigma_0,S^R(M)\otimes E\bigr)\rightarrow L^2\bigl(\Sigma_1,S^R(M)\otimes E\bigr)$ is the evolution operator associated to $D$, i.e., the unitary operator that maps $f|_{\Sigma_0}$ to $f|_{\Sigma_1}$ for $f\in\Ker(D)$.
			\item $A_i$ is the twisted Dirac operator over the Cauchy hypersurface $\Sigma_i$. We identify the twisted spinor bundle of $\Sigma_i$ with the restriction of the right-handed twisted spinor bundle of $M$ using the past directed unit normal vector field to $\Sigma_i$.
		\end{itemize}
	\end{Setting}
	We can equip the Dirac operator with various boundary conditions.
	\begin{Definition}
		For projections $P_i$ in $L^2\bigl(\Sigma_i,S^R(M)\otimes E\bigr)$, let $D_{P_0,P_1}$ be the Dirac operator with boundary conditions given by the ranges of $P_i$. This means we take (the closure of) $D$ over~$M_0$ and restrict its domain to those spinor fields $\psi$ such that $\psi|_{\Sigma_i}\in \Ran(P_i)$.
	\end{Definition}
	 An important set of boundary conditions are the so-called $APS$ boundary conditions, given by the projections onto the negative/positive eigenspaces of $A_0$ (resp.\ $A_1$).
	\begin{Definition}
		Define the following spectral projections (with $\chi_I$ denoting the characteristic function of an interval $I$):
		\[P_-(0):=\chi_{(-\infty,0)}(A_0)\qquad \text{and}\qquad P_+(1):=\chi_{(0,\infty)}(A_1).\]
		The boundary conditions defined by these projections are known as Atiyah--Patodi--Singer (APS) boundary conditions and we also write $D_{APS}:=D_{P_-(0),P_+(1)}$.
	\end{Definition}
	\begin{Remark}
		We can also allow $P_-(0)$ and $P_+(1)$ to be any finite rank perturbation of these projections. This will not affect any of the properties needed for our purposes. In particular, we could replace them by $\chi_{(-\infty,a)}(A_0)$ and $\chi_{(b,\infty)}(A_1)$ for any $a,b\in \mathbb{R}$ and still obtain the same results. Moreover, we could also interchange positive and negative spectral projections. Any change that preserves the statement of Theorem~\ref{APScomp} will not affect any of the further considerations.
	\end{Remark}

	We shall use two crucial results about the Lorentzian setting, that are also essential in~\cite{BH}. The rest will only be abstract considerations about Fredholm pairs.
	Firstly, the index of $D_{P_0,P_1}$ can be expressed in terms of the boundary projections and the evolution operator.
	\begin{Theorem}
		\label{bdyindex}
		The operator $D_{P_0,P_1}$ is Fredholm with index $k$ if and only if
        $
        \bigl(Q^{-1}P_1Q,1-P_0\bigr)$
        is a Fredholm pair of index~$k$.
	\end{Theorem}
	This is essentially \cite[Theorem 4.2]{BH}, using that the index of the pair of subspaces $(\Ran(P_0),\allowbreak Q\Ran(P_1))$ as defined there is exactly the index of $\bigl(Q^{-1}P_1Q,1-P_0\bigr)$.
	Roughly speaking, this theorem means the Dirac operator is Fredholm if and only if the boundary projections are ``almost complementary'' when related via time evolution.

	Secondly, we need that the APS boundary conditions are ``as Fredholm as they can be'', i.e., the relevant projections have compact difference.
	\begin{Theorem}	\label{APScomp}
		The difference
$1-P_-(0)-Q^{-1}P_+(1)Q$
		is a compact operator.
	\end{Theorem}
	\begin{Remark}
		This means in particular that $\bigl(Q^{-1}P_+(1)Q,1-P_-(0)\bigr)$ and hence $D_{P_-(0),P_+(1)}$ is Fredholm.
	\end{Remark}
	\begin{proof}[Proof of Theorem~\ref{APScomp}]
		This is a fairly easy reformulation of \cite[Lemma 2.6]{BSglobal}, which states that
		\[Q_{+-}:=P_+(1)QP_-(0)\qquad \text{and}\qquad Q_{-+}:=(1-P_+(1))Q(1-P_-(0))\]
		are compact operators.
		We have
		\begin{align*}
			1-P_-(0)-Q^{-1}P_+(1)Q={}&Q^{-1}(Q(1-P_-(0))-P_+(1)Q)\\
            ={}&Q^{-1}(P_+(1)Q(1-P_-(0))+(1-P_+(1))Q(1-P_-(0))\\
			&{}{-}\,(P_+(1)QP_-(0)+P_+(1)Q(1-P_-(0)))\\
			={}&Q^{-1}(Q_{-+}-Q_{+-}).
		\end{align*}
		As the compact operators are an ideal, this is compact.
	\end{proof}

	\section{Fredholm pairs}
	\label{s3}
	We derive some further criteria for the Fredholmness of a pair, eventually leading to the characterization that we want. Throughout this section, let $P$ and $P'$ be projections in a complex Hilbert space $H$. For later use, we first compute the difference of certain two-dimensional projections.
	\begin{Lemma}
		\label{2dE}
		Let $p$ be the projection onto $\spann\{(a,b)\}$ in $\mathbb{C}^2$, for $a,b\in \mathbb{R}$, and let $q$ be the projection onto $\mathbb{C}\times0$. Then $\|p-q\|=\frac{|b|}{\sqrt{a^2+b^2}}$.
	\end{Lemma}
	\begin{Remark}
		This is the sine of the angle between the ranges of $p$ and $q$.
	\end{Remark}
	\begin{proof}[Proof of Lemma \ref{2dE}]
		Assume without loss of generality that $\|(a,b)\|^2=a^2+b^2=1$, as the result is invariant under rescaling $(a,b)$.
		The matrix representation of $p-q$ is
		\[
        p-q=(a,b)\otimes(a,b)-\begin{pmatrix}{1}&{0}\\{0}&{0}\end{pmatrix}
        =\begin{pmatrix}{a^2-1}&{ab}\\{ab}&{b^2}\end{pmatrix}
        =b\begin{pmatrix}{-b}&{a}\\{\hphantom{-} a}&{b}\end{pmatrix}.
        \]
		As the last matrix (without prefactor) is unitary, the norm of this is given by $|b|$ as desired.
	\end{proof}

	\begin{Remark}[two-dimensional model]
In order to give the reader a better intuition of what is going on, we will look at what happens for direct sums of two-dimensional projections. As~we are working with complex Hilbert spaces, we want to consider projections in $\mathbb{C}^2$, however for pictures and geometric intuition it is easier to think of projections in $\mathbb{R}^2$. We will thus consider projections of the form $p\otimes 1$ in $\mathbb{C}^2\cong \mathbb{R}^2\otimes \mathbb{C}$.

		On the Hilbert space \smash{$\bigoplus_{i=1}^\infty \mathbb{C}^2$}, we consider projections
\[
P:=\bigoplus\limits_{i=1}^\infty p_i\otimes 1
\qquad \text{and}\qquad
P':=\bigoplus\limits_{i=1}^\infty p_i'\otimes 1
\]
for $p_i$, $p_i'$ projections on a line in $\mathbb{R}^2$. The pair will be Fredholm if and only if
$
\bigl\|p_i-p_i'\bigr\|<1-\varepsilon$
for all but finitely many $i$ (as in that case $\|P-P'\|_C<1-\varepsilon$).		
		This means that any condition enforcing Fredholmness of $P-P'$ must prevent the ranges of $p_i$ and $p_i'$ from coming arbitrarily close to being mutually orthogonal for all but finitely many $i$, i.e., the angle  between these ranges should be bounded by a constant less than \smash{$\frac{\pi}{2}$}.
	\end{Remark}
\medskip
	For arbitrary projections, the condition ``for all but finitely many $i$'' in the remark translates into finite codimension subspaces.
	
\begin{Lemma}\label{11}
		The following are equivalent:
		\begin{enumerate}\itemsep=0pt
\item[$(1)$] The pair $({P'},P)$ is a Fredholm pair.
\item[$(2)$] There is $\varepsilon>0$ and a closed subspace $V$ of finite codimension such that  ${\|{P'}x\| \geq \varepsilon\|x\|}$ for $x\in \Ran(P)\cap V$ and $\|Px\|\geq\varepsilon\|x\|$ for $x\in \Ran({P'})\cap V$.
\item[$(3)$] There are closed finite-codimensional subspaces $V_1\subseteq \Ran(P)$ and $V_2\subseteq\Ran({P'})$ $($finite-codimensional in the respective ranges, not in $H)$ and $\varepsilon>0$ such that we have  $\|{P'}x\|\geq\varepsilon\|x\|$ for $x\in V_1$ and $\|Px\|\geq\varepsilon\|x\|$ for $x\in V_2$.
\end{enumerate}
\end{Lemma}
	\begin{proof}
		$(1) \Rightarrow (2)$: By Proposition~\ref{FredDiff}, $(P',{P})$ is Fredholm if and only if $\|P-{P'}\|_C<1$. This means that there is a closed $V$ of finite codimension such that $\|P-{P'}|_V\|<1-\varepsilon$ for some $\varepsilon\geq 0$. For $x\in V\cap\Ran(P)$, we have
		\[\bigl\|{P'}x\bigr\|\geq\|Px\|-\bigl\|\bigl(P-{P'}\bigr)x\bigr\|\geq\|x\|-(1-\varepsilon)\|x\|=\varepsilon\|x\|.\]
		For $x\in V\cap \Ran({P'})$, we proceed analogously.
		
		$(2) \Rightarrow (3)$: Choose $V_1:=V\cap \Ran(P)$ and $V_2:=V\cap \Ran({P'})$. As intersecting with a~subspace cannot increase the codimension (when taking that subspace as the new total space) and intersections of closed spaces are closed, these are closed, finite-codimensional and have the desired properties.
				
		$(3)\Rightarrow (1)$:	
		$P_{P'}$ has finite-dimensional kernel, as the kernel of $P$ cannot intersect $V_2$ non-trivially. As the kernel of ${P'}$ cannot intersect $V_1$ non-trivially, we obtain that ${P'}_P$ has finite-dimensional kernel. As this is the adjoint of $P_{P'}$, the only thing left to show is that $P_{P'}$ has closed range.
		
		Let $(z_n)$ be a sequence in $\Ran(P_{P'})$ that converges to a limit $z\in\Ran(P)$. We need to show that $z\in P(\Ran(P'))$. Assume without loss of generality that the $z_n$ are non-zero. Choose $x_n\in \Ran({P'})$ orthogonal to $\Ker(P_{P'})$ such that $Px_n=z_n$ (this can be obtained from an arbitrary preimage sequence by subtracting the projection to the kernel). Choose $a_n$ and $b_n$ in $\Ran({P'})$ such that $a_n\in V_2$, $b_n\bot V_2$ and $a_n+b_n=x_n$.
		
\begin{Claim}
The sequence $(x_n)$ is bounded.
\end{Claim}
		
		Assume $(x_n)$ were unbounded. Without loss of generality, assume that $\|x_n\|\rightarrow\infty$, otherwise restrict to a subsequence. $Px_n\rightarrow z$, so in particular $(Px_n)$ is bounded. Thus we have
		\[P\frac{x_n}{\|x_n\|}\rightarrow 0.\]
		As $\|a_n\|^2+\|b_n\|^2=\|x_n\|^2$ by orthogonality, we know that \smash{$\bigl(\frac{b_n}{\|x_n\|}\bigr)$} is bounded. Thus it has a convergent subsequence converging to some limit $\beta$, as the orthogonal complement of $V_2$ in $\Ran({P'})$ is finite-dimensional. By restricting to a subsequence, we may assume without loss of generality that \smash{$\bigl(\frac{b_n}{\|x_n\|}\bigr)$} converges to $\beta$. Then
\[
P\frac{a_n}{\|x_n\|}=P\frac{x_n}{\|x_n\|}-P\frac{b_n}{\|x_n\|}\rightarrow-P\beta.
\]
		In particular, $\bigl(P\frac{a_n}{\|x_n\|}\bigr)$ is Cauchy. As $P$ is bounded below on $V_2$, $\bigl(\frac{a_n}{\|x_n\|}\bigr)$ is also Cauchy and converges to some limit $\alpha$, which must satisfy $P\alpha=-P\beta$. We obtain that
\[
\frac{x_n}{\|x_n\|}\rightarrow \alpha+\beta\in \Ker(P_{P'}).
\]
		As $\frac{x_n}{\|x_n\|}$ are unit vectors orthogonal to $\Ker(P_{P'})$, this is a contradiction, which finishes the proof of the claim.
		
		We now use a similar argument without normalization to conclude that the range is closed. We know that $(x_n)$ and thus also $(b_n)$ is bounded. As $(b_n)$ is contained in a finite-dimensional space, it has a convergent subsequence. Assume without loss of generality that $(b_n)$ converges to a limit $b$. Thus
		\[Pa_n=z_n-Pb_n\rightarrow z-P b.\]
		This means that $\left(Pa_n\right)$ is Cauchy. As $P$ is bounded below on $V_2$, $(a_n)$ is Cauchy, so it converges to some limit $a$. By continuity,
		\[P(a+b)=\lim\limits_{n\rightarrow\infty} Px_n=z.\]
		Thus any $z$ in the closure of $\Ran(P_{P'})$ is already in $\Ran(P_{P'})$, i.e., $P_{P'}$ has closed range.
	\end{proof}

	Two projections of difference smaller than 1 can always be deformed into each other by a~path of projections that does not increase the distance from one endpoint. The same can also be done with the Calkin norm instead of the full operator norm.
	\begin{Lemma}
		\label{PPath}
		Let $(P_0,P_1)$ be a Fredholm pair of projections. Then there is a continuous path of projections $P(t)$ such that $P(0)-P_0$ is compact, $P(1)=P_1$ and
		\[\|P(t)-P_0\|_C\leq \|P_1-P_0\|_C.\]
	\end{Lemma}
	\begin{proof}
		We use the construction of \cite[Propositions 5.3.18 and 5.3.19]{DSW} in order to deform $P$ to a~compact perturbation of $S$. The path we use is the same as in \cite[Proposition~5.3.19]{DSW}. We~repeat the construction and check that it satisfies the desired norm bounds.
		
		Set $F(x):=1-2x$. Given two projections $P_0$ and $P_1$ with $\|P_1-P_0\|<1$, define a path by
		\begin{align*}
			W(P_0,P_1)(t):={}&F^{-1}\biggl( \left(1+ (F(P_0)F(P_1)-F(P_1)F(P_0))\sin\left(\frac{\pi}{2}t\right)\cos\left(\frac{\pi}{2}t\right)\right)^{-\frac{1}{2}} \\
			& {}\times \left(F(P_0)\cos\left(\frac{\pi}{2}t\right)+F(P_1)\sin\left(\frac{\pi}{2}t\right)\right)\biggr).
		\end{align*}
		It is shown in the proof of \cite[Proposition~5.3.18]{DSW} that this is a well-defined, continuous path of projections from $P_0$ to $P_1$ satisfying
		\[\|W(P_0,P_1)(t)-P_0\|\leq\|P_1-P_0\|.\]
		The proof works in an arbitrary $C^*$-algebra, so we can also use this in the Calkin algebra.
		
		Now let $(P_0,P_1)$ be any Fredholm pair of projections.
		Set $V_0:=\Ran(P_0)\cap\Ker(P_1)$ and $V_1:=\Ker(P_0)\cap\Ran(P_1)$ and define $V':=(V_0\oplus V_1)^\bot$. As the $V_i$ either lie in the image or the kernel of the $P_i$, the restrictions of the $P_i$ to $V'$ are projections in $V'$. Moreover, as we removed the $\pm1$-Eigenspaces of $P_1-P_0$ and $\pm1$ were not in the essential spectrum of $P_1-P_0$, we have
		\[\|P_1-P_0|_{V'}\|<1.\]
		Define
		\[P(t):= W(P_0|_{V'},P_1|_{V'})\oplus 0_{V_0}\oplus 1_{V_1}.\]
		As $W(P_0|_{V'},P_1|_{V'})$ is a continuous path of projections, so is $P(t)$. We also have $P(1)=P_1$ and $P(0)-P_0=0_{V'}\oplus 1_{V_0}\oplus -1_{V_1}$ is compact. It remains to check the norm bound.
		
		Let $\pi$ denote the projection onto the Calkin algebra. For operators on $V'$, define \[\tilde\pi(T):=\pi(T\oplus 0_{V_0\oplus V_1}).\] As $V_0$ and $V_1$ are finite-dimensional, we get the same result when extending with something other than zero. Moreover, $\tilde\pi$ is a $*$-homomorphism, so everything used to construct $W$ is preserved. We get
		\begin{align*}
			\|P(t)-P_0\|_C&=\|\pi(P(t)-P_0)\|
			 =\|\tilde\pi (W(P_0|_{V'},P_1|_{V'})(t)-P_0|_{V'})\| \\
			 &=\|W(\tilde\pi(P_0|_{V'}),\tilde\pi(P_1|_{V'}))(t)-\tilde\pi(P_0|_{V'})\|\\
			&\leq \|\tilde\pi(P_1|_{V'})-\tilde\pi(P_0|_{V'})\|
			 =\|\pi(P_1)-\pi(P_0)\|
			 =\|P_1-P_0\|_C. \tag*{\qed}
        \end{align*}
\renewcommand{\qed}{}
\end{proof}

		We want to show that a pair of projections is Fredholm if both projections are not too far away from some reference projection $S$.
	Thus we let $S$ be another projection and let $R:=P-S$ and $R':=P'-S$ denote the differences from this. In our application, $S$ will correspond to one of the APS boundary projections. A weaker version of this (without the squares) was shown in \cite[Proposition~5.3.15]{DSW}.
	\begin{Theorem}
		\label{Cnormfred}
		If $\|R\|_C^2+\bigl\|R'\bigr\|_C^2<1$,  then $(P,P')$ is Fredholm and
		\[ \ind\bigl(P,P'\bigr)=\ind(P,S)+\ind\bigl(S,P'\bigr).\]
	\end{Theorem}
	\begin{proof}
		Let $\varepsilon\in(0,1)$ be small enough such that
$\|R\|_C^2+\bigl\|R'\bigr\|_C^2<(1-\varepsilon)^2$.
		Choose a finite codimension subspace $V$ such that
		\[\|R|_V\|^2+\bigl\|R' |_V\bigr\|^2<(1-\varepsilon)^2.\]
		Let $x\in V\cap \Ran(P)$ with $\|x\|=1$. By the parallelogram identity, we have
		\[\bigl\|Rx+R'x\bigr\|^2+\bigl\|Rx-R'x\bigr\|^2=2\bigl(\|Rx\|^2+\bigl\|R'x\bigr\|^2\bigr)<2(1-\varepsilon)^2. \]
		Thus we must have either
		\[\bigl\|Rx-R'x\bigr\|^2<(1-\varepsilon)^2\qquad \text{or}\qquad
		 \bigl\|Rx+R'x\bigr\|^2<(1-\varepsilon)^2.\]
		In the first case, we have
		\[(1-\varepsilon)\|x\|>\bigl\|Rx-R'x\bigr\|=\bigl\|Px-P'x\bigr\|\geq\|Px\|-\bigl\|P'x\bigr\|=\|x\|-\bigl\|P'x\bigr\|,\]
		in the second case, we have
		\begin{align*}
\begin{split}
			(1-\varepsilon)\|x\|&>\bigl\|Rx+R'x \bigr\|
			 =\bigl\|\bigl(P+P'-2S\bigr)x\bigr\|
			 =\bigl\|(1-2S)x+P'x\bigr\|\\
			&\geq\|(1-2S)x\|-\bigl\|P'x\bigr\|
			 =\|x\|-\bigl\|P'x\bigr\|,
\end{split}
		\end{align*}
		as $1-2S$ is unitary.
		In either case, we can conclude $\|P'x\|>\varepsilon\|x\|$. The same works with $P$ and~$P'$ interchanged (for $x\in V\cap\Ran(P')$). Thus the pair is Fredholm by Lemma \ref{11}.
		
		It remains to show the second part of the claim. Let $P(t)$ be a path (as in Lemma \ref{PPath}) of projections such that $P(1)=P$, $P(0)-S$ is compact and
		\[\|P(t)-S\|_C\leq \|P-S\|_C.\]
		As $(P,S)$ is a Fredholm pair (i.e., has Calkin norm less than $1$), so is $(P(t),S)$. Furthermore, we have
		\[\|P(t)-S\|_C^2+\bigl\|P'-S\bigr\|_C^2\leq \|R\|_C^2+\bigl\|R'\bigr\|_C^2<1,\]
		hence $(P(t),P')$ is a Fredholm pair. As continuous families of Fredholm pairs have constant index and the Fredholm index is additive if one of the pairs involved has compact difference (see Proposition~\ref{Fredholmprops}), we obtain
		\[
        \ind\bigl(P,P'\bigr)=\ind\bigl(P(0),P'\bigr)=\ind(P(0),S)+\ind\bigl(S,P'\bigr)=\ind(P,S)+\ind\bigl(S,P'\bigr).
        \tag*{\qed}
        \]
\renewcommand{\qed}{}
\end{proof}
		
	We can see that Theorem~\ref{Cnormfred} is really the best estimate that we can hope for by looking at our two-dimensional model case again: 
  \begin{figure}[ht]\centering
		\begin{tikzpicture}
			
			\draw[red, thick] ( -3,0 ) -- (3, 0) node[anchor=north] {$\Ran(p_i)$};
			\draw[blue,  thick] ( -3,-2 ) -- (3, 2) node[anchor=south] {$\Ran(q_i)$};
			\draw[green,  thick] ( -2,-3 ) -- (2, 3) node[anchor=east] {$\Ran(p_i')$};
			
			\draw[black, thick] (1.5, 0) arc[start angle=0, end angle=33.69, radius=1.5] node[midway, anchor=east] {$\beta_i$};
			\draw[violet, thick] (2, 0) arc[start angle=0, end angle=56.31, radius=2] node[midway, anchor=west] {$\alpha_i$};
			\draw[purple, thick] ({1.5*cos(33.69)}, {1.5*sin(33.69)}) arc[start angle=33.69, end angle=56.31, radius=1.5];
			\node [text=purple] at (0.8, 0.9) {$\beta_i'$};
			
		\end{tikzpicture}
	\end{figure}

	\begin{Remark}[two-dimensional model, continued]
	 Assume again that $P$ and $P'$ are direct sums as before. Suppose that $S$ is of the same form, i.e., a direct sum of two-dimensional real projections $q_i$ as above.
	 Let $\alpha_i$ be the angle between the ranges of $p_i$ and $p'_i$ and let $\beta_i$ and $\beta_i'$ be the angles between the range of $q_i$ and that of $p_i$ resp.\ $p'_i$. By Lemma \ref{2dE}, the sines of these angles are the norms of the respective differences. $\alpha_i$ is at most $\beta_i+\beta'_i$. To get the condition $\alpha_i<\frac{\pi}{2}$, we need to demand that $\beta'_i<\frac{\pi}{2}-\beta_i$, i.e.,
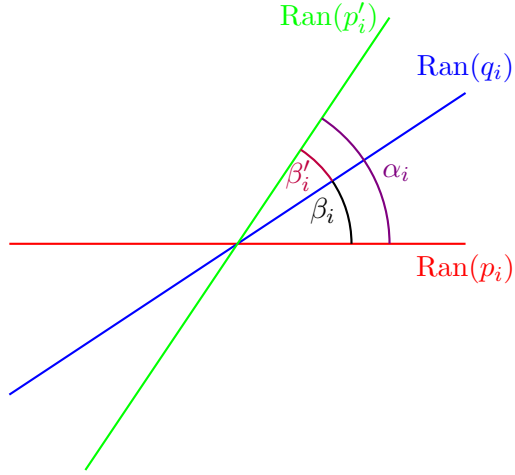
\begin{align*}
\|p_i-q_i\|^2+\bigl\|p'_i-q_i\bigr\|^2& =\sin(\beta_i)^2+\sin\bigl(\beta'_i\bigr)^2<\sin(\beta_i)^2+ \sin\left(\frac{\pi}{2}-\beta_i\right)^2\\
& =\sin(\beta_i)^2+\cos(\beta_i)^2=1.
\end{align*}
	To avoid $\alpha_i$ getting arbitrarily close to $\frac{\pi}{2}$, we need to have this uniformly, i.e.,
\[
\|p_i-q_i\|^2+\bigl\|p'_i-q_i\bigr\|^2<1-\varepsilon
\]
for a fixed $\varepsilon>0$ and all but finitely many $i$. This yields exactly the condition
\[
\|R\|_C^2+\bigl\|R'\bigr\|_C^2<1.
\]
	\end{Remark}

	We now proceed to study the case where the projections $P$ and $P'$ are graphs of linear maps with respect to a splitting of $H$ into a direct sum. These graphs can be interpreted as a~perturbation of the operators domains, with the norm of the operators providing a measure for the size of the perturbations.
	\begin{Lemma}
		\label{graphest}
		Let $H=Y\oplus Z$ be an orthogonal decomposition of $H$. Let $G\colon Y\rightarrow Z$ be a bounded linear map. Let $P$ be the projection onto the graph of~$G$, viewed as a subset of $H$, i.e.,
		\[\Gamma(G):=\{x+Gx\mid x\in Y\}.\]
		Let $S$ denote the projection onto $Y$. Then \smash{$\|P-S\|_C\leq\frac{\|G\|_C}{\sqrt{1+\|G\|_C^2}}$}.
	\end{Lemma}
	\begin{Remark}
		As the Calkin norm is unaffected by extension by $0$ and by finite rank perturbations, we may add to or remove from the domain of $G$ finite-dimensional subspaces without changing the result of this lemma and thus the subsequent theorem.
	\end{Remark}	
	\begin{proof}[Proof of Lemma \ref{graphest}]
		Before we show the estimate, there are many simplifications we can do without loss of generality.

		Firstly, we want to replace Calkin norms with operator norms. Suppose we had only shown the theorem with $\|G\|$ instead of $\|G\|_C$ on the right-hand side. Let $\varepsilon >0$ and let $Y'\subset Y$ be a~subspace of finite codimension such that
		\[\|G|_{Y'}\|\leq\|G\|_C+\varepsilon\]
		(this exists by Definition/Proposition~\ref{cnorm}). Then applying the weaker version of the lemma with~$Y'$ replacing $Y$ and $G':=G|_{Y'}$, we can deduce that
		\[\bigl\|P'-S'\bigr\|_C\leq \frac{\|G'\|}{\sqrt{1+\|G'\|^2}},\]
		where $S'$ is the projection onto $Y'$ and $P'$ that on the graph of $G'$ in $Y'\oplus Z$. Extending them by zero \big(on $Y'^\bot\subset Y$\big), we obtain the projections on these spaces in $H$, which have compact difference from $S$ (resp.\ $P$) (the ranges are finite codimension subsets of each other).
		Thus
		\[\|P-S\|_C=\bigl\|\bigl(P'-S'\bigr)\oplus 0\bigr\|_C=\bigl\|P'-S'\bigr\|_C\leq \frac{\|G'\|}{\sqrt{1+\|G'\|^2}}\leq \frac{(\|G\|_C+\varepsilon)}{\sqrt{1+(\|G\|_C+\varepsilon)^2}}.\]
		As $\varepsilon$ was arbitrary, this implies the stronger lemma. Thus it suffices to show the lemma without the subscript $C$ on the right-hand side. As we always have $\|P-S\|_C\leq \|P-S\|$, it is sufficient to show the lemma with operator norms everywhere.
		
		Secondly, we may assume that the kernel of $G$ is zero and it has dense range: As both $P$ and~$S$ are~$1$ on $\Ker(G)$ and $0$ on $\Ran(G)^\bot$, their difference is zero on either space. Thus removing these spaces does not change the left-hand side of the inequality. As it also does not change $\|G\|$ and hence the right-hand side, it suffices to show the theorem in the case where
\[
\Ker(G)=\Ran(G)^\bot=0,
\]
		which we will assume from now on.

		Thirdly, we want to replace $G$ by a self-adjoint operator. Let $G=UA$ be the polar decomposition of $G$, i.e., $A$ is self-adjoint and $U$ is an isometry from $\Ker(G)^\bot$ to $\overline{\Ran(G)}$. By the previous assumption, these spaces are everything, i.e., $U$ is a unitary map from $Y$ to $Z$. Suppose we had shown the lemma for $A$ instead of $G$ (and $Y$ instead of $Z$) and let again $P'$ and $S'$ be the projections for that case. As the graph of $G$ is the image of the graph of $A$ under $\tilde U:=1\oplus U$ (which also fixes $Y\oplus 0$), we get
		\[
        \|P-S\|=\bigl\|\tilde U\bigl(P'-S'\bigr)\tilde U^{-1}\bigr\|=\bigl\|P'-S'\bigr\|\leq \frac{\|A\|}{\sqrt{1+\|A\|^2}}=\frac{\|G\|}{\sqrt{1+\|G\|^2}}.
        \]
		Thus we may assume from now on that $Z=Y$ and $G$ is self-adjoint.

		Fourthly, we want to replace $G$ by a multiplication operator. By the spectral theorem, there is a unitary $V\colon Y\to L^2(\Omega)$ for some measure space $\Omega$ and an $L^\infty$-function $g$ such that $VGV^{-1}f=gf$ for any $f\in L^2(\Omega)$ (in particular $\|G\|=\|g\|_\infty$). As everything is invariant under conjugation by unitaries, it suffices to show the lemma for $Y=L^2(\Omega)$ and $G$ a multiplication operator.
		
		Finally, we identify $L^2(\Omega)\oplus L^2(\Omega)$ with $L^2\bigl(\Omega,\mathbb{C}^2\bigr)$.
		
		Thus, overall, have the following situation: Let $\Omega$ be a measure space with measure $\mu$ and $g\in L^\infty(\Omega)$. In $L^2\bigl(\Omega, \mathbb{C}^2\bigr)$, let $S$ be the projection onto $L^2(\Omega,\mathbb{C}\times\{0\})$ and $P$ the projection onto $\bigl\{(\psi,g\psi)\mid \psi\in L^2(\Omega)\bigr\}$. We need to show the inequality
		\[
        \|P-S\|\leq \frac{\|g\|_\infty}{\sqrt{1+\|g\|_\infty^2}}.
        \]
		Denote by $p_a$ the projection onto $\spann\{(1,a)\}$ in $\mathbb{C}^2$. Then the following define orthogonal projections onto the correct spaces:
		\begin{gather*}
        Sf(x)=p_0f(x),\qquad
        Pf(x)=p_{g(x)}f(x).
        \end{gather*}
		We can thus check for any $f\in H$, using Lemma \ref{2dE} for the pointwise estimate:
		\begin{align*}
			\|(P-S)f\|^2
			&=\int\limits_\Omega\|(p_{g(x)}-p_0)f(x)\|^2{\rm d}\mu(x)
			 \leq\int\limits_\Omega\frac{g(x)^2}{1+g(x)^2}\|f(x)\|^2{\rm d}\mu(x)\\
			&\leq\frac{\|g\|_\infty^2}{1+\|g\|_\infty^2}\int\limits_\Omega\|f(x)\|^2{\rm d}\mu(x)
			 =\frac{\|g\|_\infty^2}{1+\|g\|_\infty^2}\|f\|^2,
		\end{align*}
		which gives the desired estimate.
	\end{proof}
	
	This estimate allows us to investigate the case where both projections are given as the graph of an operator. The $S$ which we are perturbing will correspond to the APS-projections in our application.
	\begin{Remark}[two-dimensional model, continued further]
		The previous proof essentially reduces the estimates to the two-dimensional case. In this case, $\|G\|$ is the tangent of the angle~$\alpha$ between the projections, while the norm of the difference is the sine of $\alpha$.
			\end{Remark}

		Our actual use case requires one projection to form a Fredholm pair with the complement of the other. In the two-dimensional example, this transforms our question to the following: How much may we perturb two mutually orthogonal subspaces without them becoming equal?

\begin{figure}[ht]\centering
			\begin{tikzpicture}
				\draw[black, thick] ( -3,0 ) -- (3, 0) node[anchor=north] {$Y_0\approx Z_1$};
				\draw[black, thick] ( 0, -3 ) -- (0, 3) node[anchor=east] {$Y_1\approx Z_0$};
				\draw[green, thick] ( -1.5, -3 ) -- (1.5, 3) node[anchor=west] {$\Ran(P_0)$};
				\draw[green, thick] ( -1, -3 ) -- (1, 3) node[anchor=south] {$\Ran(P_1)$};
				\node [anchor=north] at (1,0) {1};
				\node [anchor=east] at (0,1) {1};
				\draw[blue] ( 1,0 ) -- (1, 2) node[midway,anchor=west] {$g_0$};
				\draw[blue] ( 0,1 ) -- (0.333,1) node[midway,anchor=south] {$g_1$};
				\draw[black, dotted] (0,2)--(1,2);
				\draw[red] (0,2) -- (0.66666,2) node[midway,anchor=south] {$g_0g_1$};
			\end{tikzpicture}
		\end{figure}
		
		If $\Ran(P_0)$ is a line through $(1,g_0)$ and $\Ran(P_1)$ is a line through $(g_1,1)$ ($g_i$ corresponding to~$\|G_i\|_C$), then the latter also contains $(g_0g_1,g_0)$. In order for $\Ran(P_1)$ to stay to the left (in the first quadrant) of $\Ran(P_0)$, we must thus have $g_0g_1<1$. This provides a more direct intuition for why we obtain the bounds we do in the following theorem.
	
\begin{Theorem}	\label {graphfred}
		Assume that $H=Y_i\oplus Z_i$ for $i=0$ and $i=1$, such that $1-S_1-S_0$ is compact, where $S_i$ denotes the projection onto $Y_i$. Let $G_i\colon Y_i\rightarrow Z_i$ be bounded operators with $\|G_0\|_C\|G_1\|_C<1$. Let $P_i$ be the projection onto the graph of~$G_i$. Then
		\begin{enumerate}\itemsep=0pt
			\item[$(1)$] $(P_1,1-P_0)$ is Fredholm.
			\item[$(2)$] $\ind(P_1,1-P_0)=\ind(S_1,1-S_0)$
		\end{enumerate}
	\end{Theorem}
	\begin{proof}
		(1) Let $a_i:=\|G_i\|_C$. By Lemma \ref{graphest}, \[\|P_i-S_i\|_C\leq\frac{a_i}{\sqrt{1+a_i^2}},\]
			so
			\[\|P_1-S_1\|_C^2+\|P_0-S_0\|_C^2\leq \frac{a_1^2}{1+a_1^2}+\frac{a_0^2}{1+a_0^2}.\]
			The condition that this be less than 1 can now be reformulated via the following equivalences:
			\begin{align*}
				&   \frac{a_1^2}{1+a_1^2}+\frac{a_0^2}{1+a_0^2}<1\\
				&\qquad{}\Leftrightarrow\ \  a_1^2\bigl(1+a_0^2\bigr)+a_0^2\bigl(1+a_1^2\bigr)<\bigl(1+a_1^2\bigr)\bigl(1+a_0^2\bigr)\\
				&\qquad{}\Leftrightarrow\ \  a_1^2+a_1^2a_0^2+a_0^2+a_0^2a_1^2<1+a_0^2+a_1^2+a_0^2a_1^2\\
				&\qquad{}\Leftrightarrow\ \  a_1^2a_0^2<1\\
				&\qquad{}\Leftrightarrow\ \ \|G_0\|_C\|G_1\|_C=a_1a_0<1.
			\end{align*}
			As the latter was true by assumption, the former inequality holds. As $(1-S_0)-S_1$ was assumed to be compact, we get
			 \begin{align*}
			 	\|P_1-S_1\|_C^2+\|(1-P_0)-S_1\|_C^2
				&=\|P_1-S_1\|_C^2+\|(1-P_0)-(1-S_0)\|_C^2\\
				&=\|P_1-S_1\|_C^2+\|P_0-S_0\|_C^2
				<1.
			\end{align*}
By Theorem~\ref{Cnormfred}, this implies that $(P_1,1-P_0)$ is Fredholm.

(2) Let $P_i(s)$ denote the projection onto the graph of $sG_i$. As $\|sG_i\|_C\leq\|G_i\|_C$ for $s\in[0,1]$, we can use the first part to conclude that $(1-P_0(s),P_1(s))$ is Fredholm. As continuous families of Fredholm pairs have constant index, we can conclude that
			\[
\ind(P_1,1-P_0)=\ind(P_1(1),1-P_0(1))=\ind(P_1(0),1-P_0(0))=\ind(S_1,1-S_0).
\tag*{\qed}
\]
\renewcommand{\qed}{}
\end{proof}

	\section{Application to Lorentzian index theory}
	\label{s4}
	With Theorems~\ref{bdyindex} and~\ref{APScomp} we can now obtain from the results of the previous section corresponding statements about the Fredholm properties of Lorentzian Dirac operators. In the setting of Section~\ref{Lorentzintro}, we have the following.

	\begin{Theorem}
		\label{finalCnorm}
		Let $P_i$ be a projection in $L^2\bigl(\Sigma_i,S^R(M)\otimes E\bigr)$ for $i\in \{0,1\}$.
		If
		\[\|P_1-P_+(1)\|_C^2+\|P_0-P_-(0)\|_C^2<1,\]
		then $D_{P_0,P_1}$ is Fredholm  and
		\[\ind(D_{P_0,P_1})=\ind(D_{APS})+\ind(P_0,P_-(0))+\ind(P_1,P_+(1)).\]
	\end{Theorem}
	\begin{proof}
		As $1-P_-(0)-Q^{-1}P_+(1)Q$ is compact by Theorem~\ref{APScomp}, we have
		\[\|P_1-P_+(1)\|_C=\bigl\|Q^{-1}P_1Q-Q^{-1}P_+(1)Q\bigr\|_C=\bigl\|Q^{-1}P_1Q-(1-P_-(0))\bigr\|_C.\]
		Moreover, we have
		\[\|P_0-P_-(0)\|_C=\|(1-P_0)-(1-P_-(0))\|_C.\]
		The condition is thus equivalent to
		\[\bigl\|Q^{-1}P_1Q-(1-P_-(0))\bigr\|_C^2+\|(1-P_0)-(1-P_-(0))\|_C^2<1.\]
		Applying Theorem~\ref{Cnormfred}, this means that
		\[\bigl(Q^{-1}P_1Q,1-P_0\bigr)\]
		is a Fredholm pair. By Theorem~\ref{bdyindex}, this means that
		$D_{P_0,P_1}$ is Fredholm.

		For the index, we calculate using Theorem~\ref{Cnormfred} and Proposition~\ref{Fredholmprops}:
		\begin{align*}
			\ind(D_{P_0,P_1})={}&\ind\bigl(Q^{-1}P_1Q,1-P_0\bigr)\\
={}&\ind\bigl(Q^{-1}P_1Q,1-P_-(0)\bigr)+\ind(1-P_-(0),1-P_0)\\
={}&\ind\bigl(Q^{-1}P_1Q,Q^{-1}P_+(1)Q\bigr)+\ind\bigl(Q^{-1}P_+(1)Q,1-P_-(0)\bigr)\\
			&{} +\ind(1-P_-(0),1-P_0) \\
={}&\ind(P_0,P_-(0))+\ind(D_{APS})+\ind(P_1,P_+(1)).
\tag*{\qed}
\end{align*}
\renewcommand{\qed}{}
\end{proof}

Using Theorem~\ref{graphfred} instead of Theorem~\ref{Cnormfred}, we obtain a generalization of \cite[Proposition~4.5]{BH}.

	\begin{Theorem}
		\label{finalgraph}
		Let $Y_0:=\Ran(P_-(0))$ and $Y_1:=\Ran(P_+(1))$. For $i\in \{0,1\}$, let
		$G_i\colon Y_i\rightarrow Y_i^\bot$
		be bounded linear maps such that
		$\|G_0\|_C\|G_1\|_C<1$.
		Let $P_i$ be the projection onto the graph of~$G_i$.
		Then $D_{P_0,P_1}$ is Fredholm with the same index as~$D_{P_-(0),P_+(1)}$.
	\end{Theorem}
	\begin{proof}
By Theorem~\ref{APScomp}, $1-P_-(0)-Q^{-1}P_+(1)Q$ is compact. $Q^{-1}P_1Q$ is the projection onto the graph of $Q^{-1}G_1Q$ \big(as an operator on $\Ran\bigl(Q^{-1}P_+(1)Q\bigr)$\big). As composing with an isometry does not change the  $\|\cdot\|_C$-norm, we have
		\[\|G_0\|_C\bigl\|Q^{-1}G_1Q\bigr\|_C<1.\]
		Thus we can apply Theorem~\ref{graphfred} with $S_0=P_-(0)$ and $S_1=Q^{-1}P_+(1)Q$ \big(i.e., with $Q^{-1}Y_1$ as the $Y_1$ in Theorem~\ref{graphfred} and $Q^{-1}G_1Q$ as $G_1$\big). This yields that $\bigl(Q^{-1}P_1Q,1-P_0\bigr)$ is Fredholm with the same index as $\bigl(Q^{-1}P_+(1)Q,1-P_-(0)\bigr)$. By Theorem~\ref{bdyindex}, this means that $D_{P_0,P_1}$ is Fredholm with the same index as $D_{P_-(0),P_+(1)}$.
	\end{proof}

\subsection*{Acknowledgements}
I would like to thank Christian B\"ar and Christoph Brinkmann for helpful discussions and feedback. Moreover, I would like to thank the anonymous referees for their detailed feedback.

\pdfbookmark[1]{References}{ref}
\LastPageEnding

\end{document}